\documentclass{amsart}

\usepackage{amssymb}

\theoremstyle{plain}

\numberwithin{equation}{section}

\newtheorem*{conjecture}{Conjecture}

\newtheorem*{teorema}{Theorem A}
\newtheorem*{teoremab}{Theorem B}

\newtheorem{lemma}[equation]{Lemma}

\newtheorem{proposition}[equation]{Proposition}

\newtheorem{claim}[equation]{Claim}

\newtheorem{comment}[equation]{}

\newcommand{\Irr}{\operatorname{Irr}}
\newcommand{\Lin}{\operatorname{Lin}}

\newcommand{\modu}{\operatorname{mod}}

\newcommand{\Ker}{\operatorname{Ker}}

\theoremstyle{definition}

\begin{document}

	\title{Products of characters and finite $p$-groups}

\author{Edith Adan-Bante}

\address{Department of Mathematics, 
University of Illinois at Urbana-Champaign,
Urbana, 61801}

\email{adanbant@math.uiuc.edu}

\keywords{Products of characters and nilpotent groups}

\thanks{This research was partially supported by grant NSF 9970030.}

\subjclass{20c}

\date{2002}

\begin{abstract}
 Let $G$ be a finite $p$-group, where $p$ is a prime number,
 and $\chi$ and  $\psi$ be  faithful complex irreducible
characters of $G$. We study the relation between the number $\eta(\chi,\psi)$
of  distinct irreducible constituents of the product
 $\chi \psi$ and the
characters $\chi$ and $\psi$.

\end{abstract}

\maketitle

\begin{section}{Introduction}

	Let $G$ be a finite $p$-group, where $p$ is a 
 prime number. Let  $\chi$ and $\psi$
 be irreducible faithful
complex characters of $G$. Since any product of characters is a character,
$\chi \psi$ is a character of $G$. Thus it can be written as
an integral linear combination of irreducible characters of $G$.
 Let $\eta(\chi,\psi)$ be the number of distinct
irreducible constituents
of the character $\chi\psi$.
 We study the relation between 
the number  $\eta(\chi,\psi)$ 
and the characters $\chi$ and $\psi$. Through this work,
we use the notation of \cite{isaacs}.

Does there  exist a $p$-group $G$, where $p=5$,
with  
faithful characters $\chi, \psi\in \Irr(G)$ 
such that the product 
$\chi\psi$ can be written as a non-trivial integral linear
 combination of $2$
distinct irreducible characters of $G$, i.e. 
$\eta(\chi, \psi)=2$? 
The answer is no, such 
groups with such characters can not exist. 
Moreover, for any prime $p\geq 5$, the answer remains no.
For a fixed prime $p$, we regard that as a ``gap'' 
among the possible values 
that $\eta(\chi, \psi)$  can take for any finite
$p$-group $G$ and any  faithful 
characters $\chi, \psi \in \Irr(G)$. More can be said in that
regard and  the main result in this work
is the following

\begin{teorema} 
Let $G$ be a finite $p$-group, where $p$ is a prime number.
Let $\chi, \psi \in \Irr(G)$ be faithful 
characters. Then either $\eta(\chi,\psi)=1$ or
$2\eta(\chi,\psi)>p$.
\end{teorema}

In section 6, we show that for every prime $p>2$ and every
integer $m>0$, there exist a $p$-group $G$ and a 
character 
$\chi\in \Irr(G)$ 
  such that $2\eta(\chi, \chi)-1=p$
and $\chi(1)=p^m$.

 We wonder if there exist other ``gaps''.
For example: does there exist a $5$-group $G$
with  
faithful characters $\chi, \psi\in \Irr(G)$ 
such that  
$\eta(\chi, \psi)=4$? More 
specifically we make the 

\begin{conjecture}
Let $G$ be a finite $p$-group, where $p$ is an odd prime number.
Let $\chi, \psi \in \Irr(G)$ be faithful 
characters. If $2\eta(\chi,\psi)-1>p$, then 
$\eta(\chi,\psi)\geq p$.
\end{conjecture}

If $\chi \in \Irr(G)$, denote by $\overline{\chi}$ its complex conjugate,
i.e. $\overline{\chi}(g) = \overline{\chi(g)}$ for all 
$g\in G$. An application of Theorem A is the following

\begin{teoremab} 
Let $n$ be a positive  integer. Then there exists a finite 
set $S_n$ of positive integers such that 
for any nilpotent group $G$ and any $\chi \in \Irr(G)$ 
with $\eta(\chi, \overline{\chi})=n$, we have
\begin{equation*}
\chi(1) \in S_n.
\end{equation*}
\end{teoremab}
   
In \cite{edith2} we proved that 
for every prime  $p$, 
 there exist a supersolvable group $G$ and 
a character $\chi \in \Irr(G)$ such that 
$\eta(\chi, \overline{\chi})=3$ and $\chi(1)= p$. Thus 
Theorem B does not hold true assuming the weaker hypothesis
that the groups are supersolvable.  

\end{section}

\begin{section}{Preliminaries}

\begin{lemma} \label{contradictionsetomega}
Let $p$ be a prime number and $\epsilon$ be a primitive $p$-th
root of unity.
Let $\Omega$ be a proper subset of $\{0,1,2, \ldots, p-1\}$ and
$\{c_i \mid i \in \Omega\}$ be a set of integer numbers.
If
\begin{equation*}
\sum_{i\in \Omega} c_i \epsilon^i =0
\end{equation*}
\noindent then $c_i =0 $ for $i \in \Omega$. 
\end{lemma}
\begin{proof} The minimum polynomial of a primitive $p$-th root
of unity  over ${\bf Q}(x)$ is 
$$p(x)=\sum_{i=0}^{p-1} x^p.$$
Set $q(x)= \sum_{i\in \Omega} c_i x^i$. Then $p(x)$ divides
$q(x)$. If $q(x)\neq 0$, there is some non-zero $r(x)\in {\bf Q}(x)$
such that $q(x)= p(x)r(x)$. 
Since the degree of $q(x)$ is at most the degree  $p-1$ of $p(x)$,
this can only happen when $r(x)=r$  is a non-zero constant. 
But then $c_i= r \neq 0$ for all $i=0, \ldots, p-1$, 
contradicting the fact that $\Omega$ is
a proper subset of $\{0,1,2, \ldots, p-1\}$.
We conclude that $q(x)=0$. Thus $c_i =0$ for all
$i$.
\end{proof}

\begin{lemma}\label{lem.1}
Let $\chi \in \Irr(G)$, where $G$ is a $p$-group.
Suppose $Z < G$, $Y$ is a normal subgroup of
$G$, that $Z< Y$ and  that $|Y:Z| = p$. If
$Z \leq {\bf Z}( \chi) $ and $Y \not\leq  {\bf Z}(\chi)$, 
then $\chi$ vanishes
on $Y \setminus Z$.
\end{lemma}
\begin{proof}
See Lemma 2.1 of \cite{us}.
\end{proof}

{\bf Notation.} We use the notation of \cite{isaacs}.
The set of irreducible complex characters of 
$G$ is denoted by $\Irr(G)$. 
We denote by $\Lin(G)$ the set of complex 
linear characters of 
$G$.
 
	Let $H$ be a subgroup of $G$ and $\lambda \in \Irr(H)$. Then 
\begin{equation*}
	\Irr(\, G\mid \lambda\,)= \{ \chi \in \Irr(G) \mid 
[\lambda, \chi_H]\neq 0\}.
\end{equation*}

Also, if $X$ and $Y$ are normal subgroups of $G$, we set
\begin{equation*}
	\Irr(X \modu Y)= \{ \gamma \in \Irr(X) \mid 
\Ker(\gamma) \geq Y \}.
\end{equation*}
\end{section}
\begin{section}{Main Lemma}
In \cite{us} we proved that if $G$ is a finite $p$-group
and  
$\chi, \psi\in \Irr(G)$ are faithful characters such that 
the product $\chi\psi$ is a multiple of an irreducible
character of $G$, then $\chi$ and $\psi$ vanish outside the 
center ${\bf Z}(G)$ of $G$. The following lemma is a 
generalization of that result.

	\begin{lemma}[Main Lemma]\label{main}
	Suppose that $p$, $G$, $\chi$, $\psi$ and $n$ satisfy

{\normalfont (3a)} $p$ is a prime,
	
{\normalfont (3b)} $ G$ is a finite $p$-group,

{\normalfont (3c)} $\chi$ and $\psi$ are irreducible characters of $G$,

{\normalfont (3d)} $\chi$ and $\psi$ are faithful, and

{\normalfont (3e)}  $2 \eta(\chi,\psi)\leq p$, where
$\eta(\chi,\psi)$ is the number of distinct irreducible constituents
of the product $\chi\psi$. 

	Then $\chi$ and $\psi $ vanish outside the center ${\bf Z}(G)$
of $G$.
\end{lemma}
\begin{proof}
Set $n=\eta(\chi,\psi)$.
Assume that the lemma is false. Then we may choose
$p$, $G$, $\chi$, $\psi$ and $n$ satisfying the hypotheses $(3\mbox{a-e})$, but
not the conclusion of the lemma, so that

\begin{comment}\label{induction}
The lemma holds true for any quintuple $p'$, $G'$, $\chi'$, $\psi'$
and $n'$ satisfying the equivalent of $(3\textup{a-e})$ with $|G'|< |G|$.
\end{comment}

We shall prove a series of claims leading to a contradiction which
will prove the lemma.

Since $\chi$ is a faithful irreducible character of $G$, we have that
${\bf Z}(G)$ is cyclic. Set ${\bf Z}(G) =Z$. Observe that 
$G\neq Z$, otherwise $\chi$ and $\psi$ are linear characters and 
therefore satisfy the conclusion of the lemma. 
Let $Y/Z$ be a chief factor of $G$.  Because $G$ is a $p$-group,
its chief factor $Y/Z$ is cyclic of order $p$ and it is 
centralized by $G$. Since $Z$ is the 
center of $G$, it follows that

\begin{comment}\label{grouph}
$H=C_G(Y)$ is a normal subgroup of index $p$ in $G$, and commutation  in 
$G$ induces a well-defined, non-singular bilinear map 
$\tau: gH, yZ \mapsto [g,y]$ of 
$G/H \times Y/Z$ into the subgroup $\Omega (Z)$
of order $p$ in $Z$.  
\end{comment}	

Let $\theta_i\in \Irr(G)$, for  $i=1,\ldots, n $, be the  
 distinct
irreducible
constituents of $\chi \psi$. Set  
\begin{equation}\label{product}
\chi \psi = \sum_{i=1}^n m_i \theta_i
\end{equation}
\noindent where $m_i>0$ is the multiplicity of $\theta_i$ in $\chi\psi$.

Since $Z$ is the center of $G$, there exist unique linear characters
$\lambda, \mu \in \Lin(Z)$ such that 
\begin{equation}\label{restriction}
\chi_Z= \chi(1) \lambda \mbox{ and } \psi_Z= \psi(1) \mu.
\end{equation}

Because $\chi$ and $\psi$ are faithful, so are $\lambda$ and $\mu$.
Furthermore \eqref{product} implies that 
\begin{equation}\label{thetai}
(\theta_i)_Z = \theta_i(1) \lambda \mu
\end{equation}
\noindent for all $i=1,\ldots,n$.
The group $Y$ has a central subgroup $Z$ with a cyclic factor group
$Y/Z$. Hence $Y$ is an abelian normal subgroup of $G$. Since 
$\lambda \in \Lin(Z)$ is faithful, it follows from \ref{grouph}
that the set $\Lin(\, Y \mid \lambda\,)$ of all extensions 
of $\lambda$ to linear characters $\rho$ of $Y$ is a single
 $G$-conjugacy class of $p= |Y:Z|$ elements. By Clifford theory
we have that 
\begin{equation}\label{restrictionchi}
\chi_Y= \frac{\chi(1)}{p} \sum_{\rho \in \Lin(Y\mid \lambda)} \rho.
\end{equation}
A similar argument shows that 

\begin{equation}\label{restrictionpsi}
\psi_Y= \frac{\psi(1)}{p} \sum_{\sigma \in \Lin(Y\mid \mu)} \sigma.
\end{equation}

\begin{claim}\label{faithful} 
Assume the notation of \eqref{product} and 
\eqref{restriction}. Then

	(i)  $\lambda \mu \in \Lin (Z)$ is faithful. 
	
	(ii) $p\geq 3$.

	(iii) For all  $i=1, \ldots, n$, we have that
$\Ker (\theta_i)=1$. Thus ${\bf Z}(\theta_i)=Z$ and 
 \begin{equation}\label{restrictionthetaiy}
(\theta_i)_Y= \frac{\theta_i(1)}{p} \sum_{\tau\in \Lin(Y\mid \lambda\mu)} 
\tau.
\end{equation}

\end{claim}

\begin{proof}
Suppose $\lambda \mu$ is not faithful. Since $Z$ is cyclic, this
implies that 
$\Omega(Z) \leq \Ker(\lambda \mu)$, where $\Omega(Z)$ is the cyclic
subgroup of $Z$ of order $p$. By \eqref{thetai} we have 
that $\Omega(Z)\leq \Ker(\theta_i)$
for $i =1,\ldots, n$. We know from \ref{grouph} that 
$[G, Y] \leq \Omega(Z) \leq \Ker(\theta_i)$. Therefore 
$Y \leq Z(\theta_i)$, and for each $i=1, \ldots , n$, we have 
$(\theta_i)_Y= \theta_i(1) \delta_i$ 
for some $\delta_i \in \Lin(Y)$. Now \eqref{product}, \eqref{restrictionchi}
and \eqref{restrictionpsi} imply that 
\begin{equation*}
\sum_{i=1}^n m_i \theta_i(1) \delta_i= \chi_Y \psi_Y
=\frac{ \chi(1)\psi(1)}{p^2} 
\sum_{\rho\in \Lin(Y\mid \lambda), \sigma \in \Lin(Y \mid \mu).}  \rho\sigma.
\end{equation*}
 
This is impossible because the leftmost expression has at most $n$ distinct
irreducible
constituents, where $2n\leq p$ by assumption, 
and the rightmost expression has precisely
$p$ distinct irreducible constituents. 
This contradiction proves that $\lambda \mu$ is a faithful 
character.

If $p=2$, then both $\lambda$ and $\mu$ restrict to the single faithful 
linear character $\epsilon$ of the group $\Omega(Z)$ of order $2$. 
Then $(\lambda \mu)_{\Omega(Z)}= \epsilon^2=1$ and 
$\lambda \mu$ is not faithful. This contradiction with the first 
statement of the claim proves the second statement.

By \eqref{product}, for $i =1,  \ldots, n$ we have 
that $\theta_i$ lies above 
$\lambda \mu$. Assume that $\Ker(\theta_j)\neq 1$ 
for some $j \in \{1, \ldots, n\}$.
Observe that $\Ker(\theta_j)$ is a normal subgroup
of $G$. Since $G$ is a nilpotent subgroup and $\Ker(\theta_j)\neq 1$,
we have
that $Z \cap \Ker(\theta_j)\neq 1$. Observe that 
$Z\cap \Ker(\theta_j) \leq \Ker(\lambda \mu)$ since 
$\theta_j$ lies above 
$\lambda \mu$ by \eqref{thetai}.
 By (i) we have that $\Ker(\lambda \mu)=1$. 
Thus for every $i \in \{1, \ldots, n\}$ we have that 
$\Ker(\theta_i)=1$.

An argument similar to that given for \eqref{restrictionchi},
shows that \eqref{restrictionthetaiy} holds.
 
\end{proof}

\begin{claim}\label{thetaiandrho} Let $\rho \in \Irr(\,Y\mid \lambda\,)$ 
and $\sigma \in \Irr(\, Y \mid \mu\,)$ be characters of $Y$.
Set  $G_{\rho\sigma}=\{ g \in G \mid (\rho\sigma)^g=\rho\sigma \}$.
Then 

(i) $\theta_i$ lies above $\rho\sigma$ for all $i=1,\ldots, n$.

(ii) $G_{\rho\sigma}= H$. 

(iii)  For each $i=1, \ldots, n$,
there exists a  unique character
 $(\theta_i)_{\rho\sigma}\in \Irr(\,H \mid  \rho\sigma\,)$
 inducing 
$\theta_i$. 
Also 
\begin{equation}\label{restricttoysigma}
((\theta_i)_{\rho\sigma})_Y= (\theta_i)_{\rho\sigma}(1) \rho\sigma.
\end{equation}

(iv) There exist unique characters $\chi_{\rho}\in \Irr(\,H \mid  \rho\,)$
and $\psi_{\sigma}\in \Irr(\, H \mid \sigma\,)$ inducing 
$\chi$ and $\psi$, respectively.

\end{claim}
\begin{proof}
Since  $\rho\sigma \in \Irr(\,Y \mid \lambda \mu\,)$,
 (i) follows from
\eqref{restrictionthetaiy}.

{\bf (ii)}
Since $H={\bf C}_G(Y)$, clearly $ H \leq G_{\rho\sigma}$.
Since $Z= Z(\theta_i)$ by Claim \ref{faithful} (iii), 
we have $ G_{\rho\sigma}< G$. 
It follows that $H= G_{\rho\sigma}$, since $|G:H|=p$ and 
$G$ is a $p$-group.

{\bf (iii)}
Since $\theta_i$ lies above $\rho\sigma$ for 
all $i=1, \ldots, n$, while $H= G_{\rho\sigma}$, by Clifford Theory 
we have that there exists a  unique character
$(\theta_i)_{\rho\sigma}\in \Irr(\,H \mid  \rho\sigma\,)$ inducing 
$\theta_i$. Since $H= G_{\rho\sigma}$ and
$(\theta_i)_{\rho\sigma}\in \Irr(\,H \mid  \rho\sigma\,)$,
we have that
 \eqref{restricttoysigma} holds

{\bf (iv)}
By \ref{grouph} we have  that $H={\bf C}_G(Y)$ is the 
stabilizer $H= G_{\rho}= G_{\sigma}$ of both $\rho$ and $\sigma$ in 
$G$.  Thus, Clifford Theory, \eqref{restrictionchi} and 
\eqref{restrictionpsi} give us unique characters
 $\chi_{\rho}\in \Irr(\,H \mid  \rho\,)$
and $\psi_{\sigma}\in \Irr(\,H \mid \sigma\,)$ inducing 
$\chi$ and $\psi$, respectively.

\end{proof}	
\begin{claim}\label{reductiontoh} 
Let $\rho \in \Irr(\,Y\mid \lambda\,)$ 
and $\sigma \in \Irr(\,Y \mid \mu\,)$ be characters of $Y$.
Assume the notation of \eqref{product} and 
Claim \ref{thetaiandrho}.
Then there exist a subset  $S$ of $\{1, \ldots , n\}$,
and integers $l_i$ such that  
\begin{equation*}
\chi_{\rho} \psi_{\sigma} = \sum_{i\in S} l_i
 (\theta_i)_{\rho\sigma},
\end{equation*}
\noindent and $0< l_i\leq m_i$ for all $i \in S$.

\end{claim}

\begin{proof}
By Clifford Theory the restrictions of $\chi_{\rho}$ and 
$\psi_{\sigma}$ to $Y$ satisfy
$(\chi_{\rho})_Y= \chi_{\rho} (1) \rho$ and  	 
$( \psi_{\sigma})_Y = \psi_{\sigma}(1) \sigma$.
Hence 
\begin{equation*}
( \chi_{\rho}\psi_{\sigma})_Y= (\chi_{\rho})_Y ( \psi_{\sigma})_Y=
\chi_{\rho} (1)\psi_{\sigma}(1) \rho\sigma.
\end{equation*}
So there are unique non-negative integers $m_{\varphi}$ such that
\begin{equation*}
\chi_{\rho} \psi_{\sigma}= \sum_{\varphi \in \Irr(\,H \mid \rho\sigma\,)}
m_{\varphi} \varphi.
\end{equation*}
Because $H$ is $G_{\rho\sigma}$, Clifford theory tells us that
the $\varphi^G$, for $\varphi \in \Irr(\,H \mid \rho\sigma\,)$, are 
the distinct members of $\Irr(\,G \mid \rho\sigma\,)$. Hence 
	
\begin{equation*}
(\chi_{\rho} \psi_{\sigma})^G= \sum_{\varphi \in \Irr(H \mid \rho\sigma)}
m_{\varphi} \varphi^G
\end{equation*}
\noindent is the unique decomposition of $(\chi_{\rho} \psi_{\sigma})^G$
as an integral linear combination of
irreducible characters $\varphi^G$ of $G$. 

The product character $\chi\psi= (\chi_{\rho})^G (\psi_{\sigma})^G$
is the sum of $(\chi_{\rho} \psi_{\sigma})^G$ and some other characters.
Since $\chi \psi = \sum_{i=1}^n m_i \theta_i$ by \eqref{product},
it follows that 
\begin{equation*}
\sum_{\varphi \in \Irr(\,H \mid \rho\sigma\,)}
m_{\varphi} \varphi^G= \sum_{i=1}^n l_i \theta_i
\end{equation*}
\noindent for some integers $l_i$ such that
$0\leq l_i\leq m_i$ for all $i=1, \ldots, n$. Set $S=\{ i \mid l_i>0 \}$.
Hence 
$m_{\varphi}=0$ for all $\varphi\in \Irr(\,H \mid \rho \sigma\,)$ except 
for those in $\{(\theta_i)_{\rho\sigma} \mid i \in S\}$, and
$m_{\varphi}=l_i$ if $\varphi=  (\theta_i)_{\rho\sigma}$. 
Thus
\begin{equation*}
	\chi_{\rho} \psi_{\sigma} = \sum_{i\in S} 
l_i (\theta_i)_{\rho\sigma}. 
\end{equation*} 
\end{proof}

If $Y=H$ then $|G: Z|=|G:H||H:Z|= p^2$, and 
$\chi(1)=\psi(1)=p$. Thus $\chi$ and $\psi$ vanish outside
$Z$ by Corollary 2.30 of \cite{isaacs}. Since $p$, 
$G$, $\chi$,
$\psi$ and  $n$ do not satisfy the conclusion of
the theorem, we have that $Y<H$.

	  If the abelian group $Y$ is not cyclic, then it
is the direct 
product $Y = Z \times K$ of its cyclic subgroup $Z$ of index $p$ with 
some subgroup 
$K$ of order $p$. Evidently $\rho= \lambda \times 1_K$ is one
of the characters $\rho$ in $\Lin (\,Y \mid \lambda\,)$ and has
$K$ as its kernel. Similarly $\sigma= \mu \times 1_K\in \Lin (\,Y \mid \mu\,)$ 
has $K$ in its kernel. Thus we may choose $\rho$ and $\sigma$
such that 
\begin{equation}\label{kernelsk}
\rho\in \Lin (Y \mid \lambda),\  \sigma \in \Lin (Y \mid \mu)
\mbox{ and } \Ker(\rho) = \Ker(\sigma) = K.
\end{equation}
If $Y$ is cyclic, then any characters
 $\rho\in \Lin (\,Y \mid \lambda\,)$, 
$\sigma\in \Lin (\,Y \mid \mu\,)$
are faithful since $\lambda$ and $\mu$ are faithful. 
So in this case we may also choose  
$\rho$ and $\sigma$ so that \eqref{kernelsk} holds with 
$K=1$.

In both cases \ref{grouph} implies that $H={\bf C}_G(Y)$ is the 
stabilizer $H= G_{\rho}= G_{\sigma}$ of both $\rho$ and $\sigma$ in 
$G$. Hence $K= \Ker(\rho)= \Ker(\sigma)$ is a normal subgroup 
of $H$. 

Now we pass to the factor group $ \overline{H}=H/K$ 
and its normal subgroup 
$\overline{Y}=Y/K$. In view of \eqref{kernelsk}, the $H$-invariant
characters
$\rho$ and $\sigma$ of $Y$ are inflated from unique faithful 
characters $\overline{\rho}$ and $\overline{\sigma}$, 
respectively, 
in $\Lin(\overline{Y})$. So $\overline{Y}$ is a cyclic central
 subgroup of $\overline{H}$. Since 
$(\chi_{\rho})_Y= \chi_{\rho} (1) \rho$ and  	 
$( \psi_{\sigma})_Y = \psi_{\sigma}(1) \sigma$, 
we have that $\chi_{\rho}$ and $\psi_{\sigma}$
are inflated from unique characters 
$\overline{\chi} \in \Irr(\, \overline{H} \mid \overline{\rho}\,)$
and $\overline{\psi} \in \Irr(\, \overline{H} \mid \overline{\sigma}\,)$.
Also $(\theta_i)_{\rho\sigma}$ inflates from 
$\overline{\theta_i}\in \Irr(\, 
\overline{H} \mid \overline{\rho}\overline{\sigma}\,)$, for each 
$i=1, \ldots , n$.
By Claim \ref{reductiontoh} we have that 
\begin{comment}\label{countingreduction}
 the number $l$ of distinct irreducible constituents of 
$\overline{\chi}\overline{\psi}$ is at most  $n$.
\end{comment}
Let $V/K = {\bf Z}(H/K)$. Note that $Y \neq G$ and
that $Y \leq {\bf Z}( H) \leq  V$.

\begin{claim}\label{notthecenter}
$V > Y$.
\end{claim}
\begin{proof}
Suppose that $\overline{Y}= {\bf Z}(\overline{H})= V/K$. 
Since $\rho$
and $\sigma$ are faithful linear characters of $\overline{Y}$,
and $G$ is a $p$-group, 
the characters
$\overline{\chi} \in \Irr(\,\overline{H} \mid \overline{\rho}\,)$
and 
$\overline{\psi} \in \Irr(\,\overline{H} \mid \overline{\sigma}\,)$
are faithful. 
So all the conditions ($3$a-e) are satisfied by
$p$, $\overline{H}$,
$\overline{\chi}$, $\overline{\psi}$ and $l$.
Since $|\overline{H}| \leq |H| < |G|$, by the inductive hypothesis
\ref{induction} we have that  
 each of the characters $\overline{\chi}$
and $\overline{\psi}$ vanishes on $\overline{H} \setminus \overline{Y}$.
It follows that $\chi_{\rho}$ and $\psi_{\sigma}$  vanish
on $H \setminus Y$. Hence 
$\chi=(\chi_{\rho})^G$ and $\psi=(\psi_{\sigma})^G$ 
vanish on $G\setminus Y$.
But \eqref{restrictionchi} and \eqref{restrictionpsi}
imply that $\chi$ and $\psi$ vanish on $Y\setminus Z$. 
This contradicts our assumption that the conclusion of the lemma
does not hold for $p$, $G$, $\chi$, $\psi$ and $n$.
Thus 
 $\overline{Y}\neq  {\bf Z}(\overline{H})$ and the 
claim is proved.
\end{proof}




\begin{claim}\label{xsubgroup}
Let $U$ be a normal subgroup of $H$ such that 
$Y<U$ and $|U:Y|=p$. Assume that 
$U \leq {\bf Z}((\theta_j)_{\rho\sigma})$  for some $j=1, \ldots, n$.
Then for all $i=1, \ldots, n$, we have that 
$U \leq {\bf Z}((\theta_i)_{\rho\sigma})$. Thus for each 
$i=1,\ldots, n$, there is some 
  $\nu_i \in \Lin (\,U \mid \rho\sigma\,)$ such that 
\begin{equation*}
 ((\theta_i)_{\rho\sigma})_U=  (\theta_i)_{\rho\sigma}(1)\nu_i.
\end{equation*}
\end{claim}
\begin{proof}
Observe that $((\theta_j)_{\rho\sigma})_U=  (\theta_j)_{\rho\sigma}(1)\nu_i$,
for some $H$-invariant 
$\nu_j \in \Lin(\,U \mid \rho\sigma\,)$, since 
 $U\leq {\bf Z}( (\theta_j)_{\rho\sigma})$.
Suppose that for some $i$ we have that 
$U \not\leq {\bf Z}((\theta_i)_{\rho\sigma})$. Since 
$Y\leq {\bf Z}(H)$, $|U:Y|=p$ and $U\trianglelefteq H$, we have that 
\begin{equation*}
((\theta_i)_{\rho\sigma})_U= \frac{(\theta_i)_{\rho\sigma}(1)}{p}
\sum_{\nu\in \Lin(\,U\mid \rho\sigma\,)}
\nu.
\end{equation*}
Thus the set $\Lin(\,U\mid \rho\sigma\,)$ forms a single $H$-conjugacy
class. In particular there is no $H$-invariant 
$\nu \in \Lin(\,U \mid \rho\sigma\,)$. This is a contradiction since 
$\nu_j$ is $H$-invariant. Thus the claim holds.  
\end{proof}
\begin{claim} 
 ${\bf Z}(H)> Y$.
\end{claim}

\begin{proof}
Certainly $Y$ is a subgroup of the center of $H= C_G(Y)$. 
We suppose that
equality occurs. Since $V > Y$, we can choose a normal subgroup $U$ 
of $H$ such that
$Y < U \leq V$ and $|U:Y| = p$. We have $1 < [H,U] \leq [H,V] \leq K$,
and thus $[H,U] = K$. In particular, we see that 
$U\leq {\bf Z}(\chi_{\rho})$ and  $U\leq {\bf Z}(\psi_{\sigma})$. Thus 
there are
unique linear characters $\alpha$ and $\beta$ in $\Lin (U)$ 
such that
\begin{equation}\label{centerchirho}
(\chi_{\rho})_U= \chi_{\rho}(1) \alpha \mbox{ and }
 (\psi_{\sigma})_U= \psi_{\sigma}(1)\beta,
\end{equation}
\noindent and so all values of $\chi_{\rho}$ and 
$\psi_{\sigma}$ on $U$ are nonzero.

For any $g\in G\setminus H$ we have that 
 $(\chi_{\rho})^g$ is an irreducible constituent of $\chi_H$.
We certainly have $K>1$, since $1< [ H,U ]\leq K$. So $|K|=p$ 
and $K \cap K^g=1$. 
Note that  $U \not\leq {\bf Z}(    (\chi_{\rho})^g )$ since 
if $U\leq {\bf Z}((\chi_{\rho})^g)$
 then $[H,U]\leq \Ker((\chi_{\rho})^g)=K^g$,
which gives the contradiction $1<[H,U]\leq K\cap K^g=1$.
Since $Y \leq {\bf Z}( H)$
and $|U:Y| = p$,  Lemma \ref{lem.1} implies that $(\chi_{\rho})^g$
 vanishes
on $U\setminus Y$. Fix $u \in U\setminus Y$. Since $(\chi_{\rho})^G= \chi$ and 
$|G:H|=p$,  we have that
\begin{equation*}
\chi(u) =  \sum_{i=1}^p (\chi_{\rho})^{g^i}(u) = \chi_{\rho} (u) \neq 0.
\end{equation*}
Similarly
\begin{equation*}
\psi(u) = \psi_{\sigma}
(u) \neq 0\, \mbox{ and }\,\theta_i(u) = (\theta_i)_{\rho\sigma}(u) \neq 0
 \mbox{ for all $i$. }
\end{equation*}

By \eqref{centerchirho} and Claim \ref{reductiontoh}, we have that 
\begin{equation}\label{14au}
((\theta_i)_{\rho\sigma})_U= (\theta_i)_{\rho\sigma}(1) \gamma, 
\end{equation}
\noindent where $\gamma=\alpha\beta\in \Lin(U)$,  for all $i\in S$.

By \eqref{product} we have
\begin{equation}\label{equation1}
\chi_{\rho}(1) \alpha(u)\psi_{\sigma}(1)\beta(u)=
\chi_{\rho} (u)\psi_{\sigma}(u) = \chi(u)\psi(u)=
\sum_{i=1}^n m_i \theta_i(u)= 
\sum_{i=1}^n m_i (\theta_i)_{\rho\sigma}(u).
\end{equation}
 By Claim \ref{reductiontoh} we have:
\begin{equation}\label{equation2}
\chi_{\rho} (1)\psi_{\sigma}(1) \alpha(u)\beta(u)=
\chi_{\rho} (u)\psi_{\sigma}(u)= \sum_{i \in S}l_i (\theta_i)_{\rho\sigma}(u)
=[\sum_{i \in S}l_i (\theta_i)_{\rho\sigma}(1)]\gamma(u).
\end{equation}

By \eqref{14au} and Claim \ref{xsubgroup}, for $i=1, \ldots, n$  we have that
$ ((\theta_i)_{\rho\sigma})_U=
(\theta_i)_{\rho\sigma}(1)  \nu_i$, 
 for some $\nu_i \in \Irr(\,U\mid \rho\sigma\,)$. Since
$(\nu_i)_Y=\rho\sigma= \gamma_Y \in \Irr(Y)$ and $|U:Y|=p$, there exists
some $\delta_i\in \Irr(U\modu Y)$ such that $\nu_i=\delta_i \gamma$.
Thus for $i=1, \ldots, n$, 
\begin{equation} \label{thetaixgamma}
 ((\theta_i)_{\rho\sigma})_U=
(\theta_i)_{\rho\sigma}(1)  \delta_i\gamma,
\end{equation} 
\noindent for some $\delta_i \in \Irr(Y\modu U)$. 
Combining this with \eqref{equation1}
and \eqref{equation2} we get 
\begin{eqnarray*}
\sum_{i=1}^n m_i (\theta_i)_{\rho\sigma}(1) \delta_i(u)\gamma(u) &= &
\sum_{i=1}^n m_i (\theta_i)_{\rho\sigma}(u)\\
&=& \chi_{\rho} (1)\psi_{\sigma}(1) \alpha(u)\beta(u)\\
&=& [\sum_{i \in S}l_i (\theta_i)_{\rho\sigma}(1)]\gamma(u).
\end{eqnarray*}
Since $\gamma(u) \neq 0$, we get
\begin{equation*}
\sum_{i=1}^n m_i (\theta_i)_{\rho\sigma}(1) \delta_i(u)=
 \sum_{i \in S}l_i (\theta_i)_{\rho\sigma}(1).
\end{equation*}
Thus
\begin{equation}\label{contradictionnormal}
\sum_{i=1}^n m_i (\theta_i)_{\rho\sigma}(1) \delta_i(u)-
\sum_{i \in S}l_i (\theta_i)_{\rho\sigma}(1)=0.
\end{equation}

Since 
 $\delta_i\in \Irr(U\modu Y)$ and $U/Y$ is cyclic of order $p$,
it follows that $\delta_i(u)$ is a $p$-th root of unit
for all  $i=1,2,  \ldots,n$. Let $\epsilon$ be a primitive $p$-th root of
unit.
Set  $\Omega_i = \{j \mid \delta_j(u) = \epsilon^i\}$
and $\Omega=\{i \mid |\Omega_i|\neq 0\}$.
Let
\begin{equation*} 
c_i= \left\{ \begin{array}{ll}
\sum_{j\in \Omega_0} m_j (\theta_j)_{\rho\sigma}(1)- 
\sum_{i \in S}l_i (\theta_i)_{\rho\sigma}(1) & \mbox{ if $i=0$}\\
\sum_{j\in \Omega_i} m_j (\theta_j)_{\rho\sigma}(1) & \mbox{ otherwise.}
\end{array}
\right.
\end{equation*}

Observe that the 
set $\{\,c_i\mid i \in \Omega\,\}$ is a set of integer numbers.
Observe that  $n<p$ implies that
 $|\{\delta_j(u)\mid j= 1, \ldots, n\}|=|\Omega|<p$. So 
$\Omega$ is a proper subset of $\{0, 1, \ldots, p-1\}$.
Thus by Claim \ref{contradictionsetomega} 
we have that
$c_i=  0$ for $i \in \Omega$. Suppose that $\Omega_i$ is non-empty
for some $i\neq 0$. Since  $m_j$ and $\theta_j(1)$ are nonzero
positive integers, for all $j=1, \ldots, n$, we have that 
 $\sum_{j\in \Omega_i} m_j (\theta_j)_{\rho\sigma}(1)\geq |\Omega_i|>0$.
 But 
$c_i=0$. Hence 
$\Omega_i$ is empty for all
$i \neq 0$. In particular  $\Omega_0 = \{0,1,\ldots, n\}$. Also, since
$c_0=0$, we have 
\begin{equation*}
\sum_{i =1}^n  m_i (\theta_i)_{\rho\sigma}(1)=
\sum_{i \in S}l_i (\theta_i)_{\rho\sigma}(1).
\end{equation*}
Thus
\begin{eqnarray*}
\chi_{\rho}(1)\psi_{\sigma}(1)
&=& \sum_{i\in S} l_i (\theta_i)_{\rho\sigma}(1)
=\sum_{i=1}^n m_i (\theta_i)_{\rho\sigma}(1)\\
&=&\frac{1}{p} \sum_{i=1}^n m_i (\theta_i)(1)
= \frac{1}{p} \chi(1)\psi(1) \\              
&=& \frac{1}{p}( p\chi_{\rho}(1))(p \psi_{\sigma}(1))
= p \chi_{\rho}(1)\psi_{\sigma}(1).
\end{eqnarray*}
Thus $p=1$, which is obviously impossible.
\end{proof}
\begin{comment}\label{xcentralh}
Fix a normal subgroup  $X$ of $G$ such that $X/Y$ is a chief factor
of $G$ and $X\leq {\bf Z}(H)$. Let
 $\alpha$, $\beta$, $\gamma \in \Lin (X)$, where
$\gamma=\alpha\beta$,
be the unique linear characters  such that
\begin{equation*}
(\chi_{\rho})_X= \chi_{\rho}(1) \alpha, 
 (\psi_{\sigma})_X= \psi_{\sigma}(1)\beta, \mbox{ and } 
((\theta_i)_{\rho\sigma})_X= (\theta_i)_{\rho\sigma}(1) \gamma
\mbox{ for } i\in S. 
\end{equation*}  
Also, for each $i=1, \ldots, n$, there is
some $\nu_i \in \Lin(\, X \mid \rho\sigma\,)$ such that
\begin{equation*}
((\theta_i)_{\rho\sigma})_X = (\theta_i)_{\rho\sigma}(1)\nu_i.
\end{equation*}
\end{comment} 
\begin{claim}\label{xgandy}
The subgroup $[X,G]$ generates $Y = [X, G]Z$ modulo $Z$. 
\end{claim} 
\begin{proof}
Since $ Y$, $X$ are normal subgroups of $G$ with 
$Y\triangleleft X$ and $|X/Y|=p$, the chief factors $X/Y$ of the 
$p$-group $G$ is centralized by $G$. 
So $[X,G]\leq Y$. Suppose that $[X,G]Z< Y$.
Since $|Y/Z|=p$, we must have       
$[X,G]\leq Z= {\bf Z}(G)$.
So commutation  in $G$ induces a bilinear map 
$$d: x Z, g {\bf C}_G(X) \mapsto [x,g]$$
\noindent of $(X/Z)\times (G/ {\bf C}_G(X))$ into 
the cyclic group $Z$. This map $d$ is non-singular
on the right by the definition of ${\bf C}_G(X)$. It is non-singular
on the left since $Z={\bf Z}(G)$. 
 Because $\alpha \in \Lin( X \mid \rho)$
extends the faithful character $\lambda\in \Irr(Z)$, and $|X:Z|=p^2$,
this implies that ${\bf C}_G(X)=G_{\alpha}$ has index $p^2$ in 
$G$. But $H$ fixes $\alpha$ by \ref{xcentralh}, since $\chi_{\rho}\in \Irr(H)$.
 Therefore $\alpha$ has, at the same 
time, $p^2$ distinct $G$ conjugates, and at most 
$p=|G:H|$ such conjugates. This contradiction proves the claim.
\end{proof}

\begin{comment}\label{choosingcoset}
Observe that  $G/H$ is cyclic of order $p$. So we may choose $g \in G$
such that the distinct cosets of $H$ in $G$ are $H$, $Hg$, $Hg^2$, $\ldots$,
$Hg^{p-1}$.
\end{comment}

 Since $\chi= (\chi_{\rho})^G$ is induced from 
$\chi_{\rho}\in \Irr (H)$, it follows from \ref{xcentralh} that
\begin{equation*}
\chi_X= \chi_{\rho}(1)( \alpha + \alpha^g + \ldots + \alpha^{g^{p-1}})=
\chi_{\rho}(1)\sum_{i=0}^{p-1} \alpha^{g^i} .
\end{equation*}
Similarly, we have that 
\begin{equation*}
\psi_X= \psi_{\sigma}(1)( \beta + \beta^g + \ldots + \beta^{g^{p-1}})
=   \psi_{\sigma}(1) \sum_{j=0}^{p-1} \beta^{g^j}.
\end {equation*}
Combining the two previous equations we have that
\begin{equation}\label{intermedia}
\chi_X \psi_X=  (\chi_{\rho}(1)\sum_{j=0}^{p-1} \alpha^{g^j})(
 \psi_{\sigma}(1) \sum_{j=0}^{p-1} \beta^{g^j})= 
 \chi_{\rho}(1) \psi_{\sigma}(1)
\sum_{i=0}^{p-1}\sum_{j=0}^{p-1} \alpha^{g^i} \beta^{g^j}. 
\end{equation}
By \eqref{product} we have that
\begin{equation}\label{comparing}
 (\chi\psi)_X = (\sum_{i=1}^n m_i \theta_i)_X=
\sum_{i=1}^n m_i[(\theta_i)_{\rho\sigma}(1)
\delta_i\sum_{j=0}^{p-1}(\alpha\beta)^{g^j}].
\end{equation}

\begin{claim}\label{orbits} Let $g\in G$ be 
 as in \ref{choosingcoset}.
For each $i=0,1,\dots ,p-1$, there exist $j \in \{0,1,\ldots, p-1\}$
and $\delta_{g^i} \in \Irr( X \modu Y)$  such that 
\begin{equation}\label{productitermsdelta}
\alpha \beta^{g^i} = (\alpha\beta)^{g^j} \delta_{g^i}.
\end{equation}
Also $|\{\delta_{g^i} \mid i=0,1,2,\ldots, p-1\}| \leq n$.
\end{claim}
\begin{proof} 
Since $H$ is a normal subgroup of $G$ 
and $\theta_i= ((\theta_i)_{\rho\sigma})^G$,
where $(\theta_i)_{\rho\sigma}\in \Irr(H \mid \rho\sigma)$,
 we have 
\begin{equation}\label{thetarestrictedtoh}
(\theta_i)_H= \sum_{j=0}^{p-1}((\theta_i)_{\rho\sigma})^{g^j}.
\end{equation}
Since  $X/Y$ is cyclic, by \ref{xcentralh} 
it follows that there exists some $\delta_i \in \Irr( X \modu Y)$
such that $\nu_i= \alpha \beta\delta_i$.
Since $X/Y$ is a chief factor of $G$,
 we have that 
$G$ acts trivially on $X/Y$. Thus
\begin{equation}\label{nu}
	 (\nu_i)^{g^j} = (\alpha\beta)^{g^j} \delta_i
\end{equation}
\noindent for all 
$j =0, 1, \ldots, p-1$. Therefore
\begin{alignat*}{2}
	(\theta_i)_X & =  ((\theta_i)_H)_X = 
(\sum_{j=0}^{p-1}((\theta_i)_{\rho\sigma})^{g^j})_X && \mbox{
 by \eqref{thetarestrictedtoh} }\\
& =  \sum_{j=0}^{p-1}(((\theta_i)_{\rho\sigma})^{g^j})_X 
 =   \sum_{j=0}^{p-1}(((\theta_i)_{\rho\sigma})_X)^{g^j} && \mbox{ 
since $X\unlhd G$ }\\
& =  \sum_{j=0}^{p-1}((\theta_i)_{\rho\sigma}(1)\nu_i)^{g^j}
 =  (\theta_i)_{\rho\sigma}(1) \sum_{j=0}^{p-1} (\nu_i)^{g^j}&& \\
& =  (\theta_i)_{\rho\sigma}(1) \sum_{j=0}^{p-1} 
( \alpha\beta\delta_i)^{g^j}=(\theta_i)_{\rho\sigma}(1) 
\sum_{j=0}^{p-1} (\alpha\beta)^{g^j}\delta_i
&& \mbox{ by \eqref{nu} } \\
& = (\theta_i)_{\rho\sigma}(1)\delta_i\sum_{j=0}^{p-1}(\alpha\beta)^{g^j}. &&
\end{alignat*}

Since $X$ is an abelian group, its irreducible characters are linear,
and any product of linear characters is irreducible. 
Observe that 
 $\chi_X \psi_X = (\chi\psi)_X$. Thus,
given $i=0,1,\ldots, p-1$, 
by \eqref{intermedia} and \eqref{comparing}
there exist some $j$ and some $k$ such that
$$\alpha\beta^{g^i}= (\alpha\beta)^{g^j}\delta_k.$$
Thus \eqref{productitermsdelta} holds with 
$\delta_{g^i}=\delta_k\in \Irr(X \modu Y)$.
Also, the set $\{\delta_{g^i} \mid i=0,1, \ldots, p-1\}$ is a subset
of 
$\{\delta_i \mid i =1, \ldots, n\}$. Therefore the last statement of 
Claim \ref{orbits} holds.
\end{proof}

\begin{claim}\label{choosingright} 
Let $g \in G$ be
 as in \ref{choosingcoset}.
Then there exist three distinct integers 
 $i,j, k \in \{0,1,2,\ldots, p-1 \}$,
and some $\delta \in \Irr(X \modu Y)$, such that 
\begin{eqnarray*}
\alpha \beta^{g^i} & = & (\alpha\beta)^{g^r} \delta, \\ 
\alpha \beta^{g^j} & = &  (\alpha\beta)^{g^s} \delta, \\
\alpha \beta^{g^k} & = &  (\alpha\beta)^{g^t} \delta,
\end{eqnarray*}

\noindent for some $r,s,t \in \{0,1,2,\ldots, p-1 \}$.
\end{claim}
\begin{proof}
By Claim \ref{orbits}, for each  $i=0,1,2, \ldots, p-1$
there exists $\delta_{g^i}\in \Irr( X \modu Y)$ such that 
$$\alpha \beta^{g^i}  = (\alpha\beta)^{g^r}\delta_{g^i},$$
\noindent for some $r \in \{0,1,\ldots, p-1\}$.
Since $\{\delta_{g^i}\mid i=0,1, \ldots, p-1\}$ has at most $n$ elements
and $n\leq \frac{p-1}{2}$, there must exist three distinct $i$,$j$,$k$
$\in \{0,1,\ldots, p-1\}$  such that $\delta_{g^i}=\delta_{g^j}=\delta_{g^k}$.
\end{proof}

\begin{claim}\label{choosingcosetgood}
We can choose the element  $g$ in \ref{choosingcoset}
such that one of the following 
holds:

(i) There exists some $j=2, \ldots, p-1$ such that
\begin{eqnarray*}
\alpha \beta^{g} & = & (\alpha\beta)^{g^r},\\
\alpha \beta^{g^j} & = &  (\alpha\beta)^{g^s},
\end{eqnarray*}
\noindent for some $r,s \in \{0,1,\ldots, p-1\}$ with $r \neq 1$.

(ii) There exist $j$ and $k$ such that $1<j<k<p$, and 
\begin{eqnarray*}\label{1,j,k}
\alpha \beta^{g} & = & (\alpha\beta)^{g^r}\delta,\\ 
\alpha \beta^{g^j} & = & (\alpha\beta)^{g^s}\delta, \\
\alpha \beta^{g^k} & = &  (\alpha\beta)^{g^t}\delta,
\end{eqnarray*}
\noindent for some  $\delta \in \Irr( X \modu Y)$ and some 
$r,s,t\in \{0,1,\ldots, p-1\}$ with $r\neq 1$.

\end{claim}
\begin{proof} Since $\alpha_Y=\rho\in \Lin(Y)$ 
and $\rho^f \neq \rho$ for any 
$f\in G\setminus H$, we have that 
$(\alpha \beta^f)_Y \neq ( \alpha^f \beta^f)_Y$. 
So if $\alpha \beta^{f}=(\alpha\beta)^{f^r}\delta$ for 
some $\delta\in \Irr(X \modu Y)$, some $f\in G\setminus H$ and some 
integer $r$, 
then $r \neq 1$.

Assume that $\alpha\beta= (\alpha\beta)^f\delta$, for some 
$f\in G$ and $\delta\in \Irr(X \modu Y)$. Since $\alpha_Y= \rho$ and 
$\beta_Y = \sigma$, we have that
$$\rho\sigma= (\alpha\beta)_Y =  ((\alpha\beta)^f\delta)_Y = 
(\rho\sigma)^f.$$            
 \noindent Since
 $G_{\rho\sigma}=H$ 
by Claim \ref{thetaiandrho} (ii), we have that $f \in H$. 
Therefore $\delta =1$.

If in Claim \ref{choosingright} we have that 
 $0 \in \{i,j,k\}$, then $\delta=1$. Assume that $i=0$.
 Since
$0< j <  p$, we have that $g^j \in G\setminus H$. Set $f= g^j$. 
Then for some $k', s', t'\in \{0,1,2, \ldots ,p-1\}$ we have that
${f}^{k'} H= g^k H$,  
${f}^{s'} H= g^s H$,
 and ${f}^{t'}H= g^t H$. Thus we have 
\begin{eqnarray*}
\alpha \beta^{f} & = & (\alpha\beta)^{{f}^{s'}}\\
\alpha \beta^{{f}^{k'}} & = &  (\alpha\beta)^{{f}^{t'}}.
\end{eqnarray*}
\noindent We conclude that if $i =0$
in Claim \ref{choosingright}, then  Claim \ref{choosingcosetgood} (i)
holds.
 
If $0 \not\in \{i,j, k\}$, set $f=g^i$.
 An  argument similar to that given for (i) 
shows that 
 Claim \ref{choosingcosetgood} (ii) holds.
\end{proof}

Let $g$ be as in Claim \ref{choosingcosetgood}.
Since $X/Y$ is cyclic of order $p$, we may choose $x \in X$ such 
that $X = Y<x>$.  By  \ref{xcentralh} we have $[X,H]=1$. 
Suppose that $[x,g^{-1}]\in Z$. Then $x$ centralizes both 
$g^{-1}$ and $H$ modulo $Z$. Hence $xZ \in {\bf Z}(G/Z)$,
which is false by Claim \ref{xgandy}. Hence 
$[x, g^{-1}] \in Y\setminus Z$ and so  
\begin{equation}\label{18}
Y=Z<y> \mbox{ is generated over } Z \mbox{ by } y=[x,g^{-1}].
\end{equation}
Since $[Y,G]\leq Z$ we have that $z=[y,g^{-1}]  \in Z$. If
$z=1$, then $G=H<g>$ centralizes $Y = Z<y>$, since $H$ centralizes 
$Y<X$ by \ref{xcentralh}, and $G$ centralizes Z. This
is impossible because $Z={\bf Z}(G)< Y$. Thus 
\begin{equation}\label{19}
z= [y,g^{-1}] \mbox{ is a non-trivial element of }Z.
\end{equation}
By \eqref{18} we have  $y=[x,g^{-1}]=x^{-1} x^{g^{-1}}$. By
\eqref{19} we have $z= [y, g^{-1}]=y^{-1} y^{g^{-1}}$ . Finally
$z^{g^{-1}}= z$ since 
$z\in Z$. Since $X=Z<x,y>$ is abelian by \ref{xcentralh},
it follows that 
\begin{equation}\label{21}
z^{g^{-j}}=z, \  y^{g^{-j}}= yz^j \mbox{ and }
 x^{g^{-j}}=xy^j z^{\binom{j}{2}},
\end{equation}
\noindent for any integer $j= 0, 1, \ldots, p-1$. 
Because $g^{-p} \in H$ centralizes $X$ by \ref{xcentralh}, 
we have
\begin{equation*} 
z^p=1 \mbox{ and } y^p z^{\binom{p}{2}}=1.
\end{equation*}
But $p\geq3$ is odd by Claim \ref{faithful}. Hence $p$ divides
$\binom{p}{2}= \frac{p(p-1)}{2}$ and $z^{\binom{p}{2}}=1$. Therefore

\begin{equation}\label{22}
y^p= z^p =1.
\end{equation}

It follows that $y^i$, $z^i$ and $z^{\binom{i}{2}}$ 
depend only on the residue of $i$ modulo 
$p$, for any integer $i\geq 0$.

\begin{claim}\label{contradiction}
 Suppose that
\begin{equation}\label{23}
\alpha \beta^{g} = (\alpha\beta)^{g^r}\delta, 
\end{equation}
\noindent and 
\begin{equation}\label{26}
\alpha \beta^{g^j}  = (\alpha\beta)^{g^s} \delta,
\end{equation}
\noindent  for some 
$j \in \{0, 1, \ldots, p-1\}$, $j \neq 1 $, 
 some $\delta \in \Irr(X \modu Y)$, and some 
$r, s\in \{0,1, \ldots , p-1\}$.
Then           
\begin{equation}\label{contradiction2}
\delta(x)= \beta(z)^{h(r-1)},
\end{equation}
\noindent where $2h\equiv 1 \modu p$. 
\end{claim}
\begin{proof}

Evaluating both sides of the equation \eqref{23} at $y$, we obtain
\begin{alignat*}{2}
\alpha(y)\beta(y)\beta(z) & =  
\alpha(y) \beta(y z )
 =   \alpha(y) \beta(y^{g^{-1}}) &&  \mbox{ by \eqref{21} }      \\
	& =  (\alpha\beta^{g})(y)= 
 \alpha^{g^r}(y)\beta^{g^r}(y)\delta(y) && \\
 & =  
\alpha^{g^r}(y)\beta^{g^r}(y) && 
 \mbox{ since $\delta \in \Irr(X \modu Y)$ } \\
& =  \alpha(y^{g^{-r}}) \beta(y^{g^{-r}})&& \\
& = \alpha(yz^r)\beta(yz^r) && \\
& =\alpha(y)\beta(y)\alpha(z)^r \beta(z)^r. && 
\end{alignat*}

Cancelling $\alpha(y)\beta(y) \neq 0$, we get
\begin{equation}\label{24}
\beta(z)= \alpha(z)^r \beta(z)^r.
\end{equation}

We may also evaluate both sides of \eqref{23} at $x$, obtaining 
\begin{alignat*}{2}
\alpha(x) \beta(x) \beta(y) & = \alpha(x) \beta(xy) 
=\alpha(x) \beta(x^{g^{-1}})&&  \mbox{ by \eqref{21}} \\
	& =
(\alpha \beta^g)(x)= ( (\alpha\beta)^{g^r}\delta)(x) && \\   
& =   \alpha^{g^r} (x) \beta^{g^r} (x)   \delta(x)    
	 = 
 \alpha(x^{g^{-r}})\beta(x^{g^{-r}}) \delta(x)&& \\
&=  \alpha(xy^r z^{\binom{r}{2}})
	\beta(xy^r z^{\binom{r}{2}})\delta(x)&&  \\
&=  \delta(x)  \alpha(x)\beta(x) 
\alpha(y)^r \alpha(z)^{\binom{r}{2}}\beta(y)^r \beta(z)^{\binom{r}{2}}.&&
\end{alignat*}

Cancelling $\alpha(x) \beta(x)\neq 0$ we get
\begin{equation}\label{25}
 \beta(y) = \delta(x) \alpha(y)^r \beta(y)^r \alpha(z)^{\binom{r}{2}}
 \beta(z)^{\binom{r}{2}}.
\end{equation}
Applying \eqref{26} at $y$ we get
\begin{alignat*}{2}
\alpha(y)\beta(y) \beta(z)^j
&= \alpha(y) \beta(yz^j)= \alpha(y) \beta( y^{g^{-j}})&&  \\
&= (\alpha\beta^{g^j})(y) = 
(\alpha\beta)^{g^s}(y)\delta(y) &&  \\
&= \alpha^{g^s}\beta^{g^s}(y) &&
\mbox{ since $\delta\in \Irr(X\modu Y)$} \\
& =\alpha(y^{g^{-s}}) \beta(y^{g^{-s}})= \alpha(yz^s)\beta(yz^s)&&\\
& = \alpha(y) \beta(y) \alpha(z)^s \beta(z)^s.&&
\end{alignat*}

Cancelling $\alpha(y) \beta(y) \neq 0$
we obtain
\begin{equation*}
\beta(z)^j=\alpha(z)^s \beta(z)^s= (\alpha(z)\beta(z))^s
\end{equation*}
But \eqref{24} gives
\begin{equation*}
\beta(z)^j= \alpha(z)^{jr}\beta(z)^{jr}= (\alpha(z)\beta(z))^{jr}.
\end{equation*}
The restriction $(\alpha\beta)_Z=\lambda \mu$ is faithful by 
Claim 
\ref{faithful}. The element $z \in Z$ is non-trivial by \eqref{19}
and has order $p$ by \eqref{22}. So $\alpha(z) \beta(z)=(\lambda \mu)(z)$
is a primitive $p$-th root of unity. 
Hence the equation 
$(\alpha(z) \beta(z))^s = (\alpha(z)\beta(z))^{jr}$ implies that 
\begin{equation}\label{27}
s \equiv jr \mbox{ mod } p.
\end{equation} 
We may also apply \eqref{26} at $x$ to get
\begin{eqnarray*}	
\alpha(x)\beta(x)\beta(y)^j\beta(z)^{\binom{j}{2}}&
=& \alpha(x)\beta(xy^j z^{\binom{j}{2}})=
\alpha(x) \beta(x^{g^{-j}}) \\
& = & (\alpha\beta^{g^j})(x) = 
\delta(x)\alpha^{g^s}(x) \beta^{g^s} (x) \\
&= &\delta(x)\alpha(x^{g^{-s}})\beta(x^{g^{-s}})= 
\delta(x)\alpha(xy^s z^{\binom{s}{2}})\beta (xy^s z^{\binom{s}{2}})\\
&=&\delta(x)\alpha(x)\beta(x)\alpha(y)^s\beta(y)^s \alpha(z)^{\binom{s}{2}} 
\beta(z)^{\binom{s}{2}}.
\end{eqnarray*}

Cancelling $\alpha(x)\beta(x) \neq 0$ we obtain
\begin{equation}\label{cancelling}
\beta(y)^j \beta(z)^{\binom{j}{2}}=
 \delta(x)\alpha(y)^s \beta(y)^s \alpha(z)^{\binom{s}{2}}
\beta(z)^{\binom{s}{2}}.
\end{equation}

But \eqref{25} gives
\begin{equation*}
\beta(y)^j= 
\delta(x)^j \alpha(y)^{jr} \beta(y)^{jr} \alpha(z)^{j\binom{r}{2}} 
\beta(z)^{j\binom{r}{2}}= \delta(x)^j \alpha(y)^{s} \beta(y)^{s} 
\alpha(z)^{j\binom{r}{2}} 
\beta(z)^{j\binom{r}{2}},
\end{equation*}
\noindent where the last equality holds since $s \equiv jr$ mod $p$.

Combining the previous equation with \eqref{cancelling}, we have
\begin{equation*}
 \delta(x)\alpha(y)^s \beta(y)^s \alpha(z)^{\binom{s}{2}}
\beta(z)^{\binom{s}{2}} \beta(z)^{-\binom{j}{2}}= 
\delta(x)^j \alpha(y)^{s} \beta(y)^{s} 
\alpha(z)^{j\binom{r}{2}} 
\beta(z)^{j\binom{r}{2}}
\end{equation*}

Cancelling terms and simplifying, we get
\begin{equation}\label{firststep}
\delta(x)^{j-1}= (\alpha(z)\beta(z))^{\binom{s}{2}- j\binom{r}{2}}
\beta(z)^{-\binom{j}{2}}.
\end{equation}

Since $p\geq3$ is odd and $s \equiv jr$ mod $p$, we have
$\binom{s}{2}\equiv  \binom{jr}{2}$ mod $p$. This and \eqref{22} 
 imply that
$\alpha(z)^{\binom{s}{2}} = \alpha(z)^{\binom{jr}{2}}$ and 
$\beta(z)^{\binom{s}{2}}= \beta(z)^{\binom{jr}{2}}$.
Hence
\begin{equation}\label{jrequation}
(\alpha(z)\beta(z))^{\binom{s}{2}- j\binom{r}{2}}=
(\alpha(z)\beta(z))^{ \binom{jr}{2}-j\binom{r}{2}}.
\end{equation}
We can check that
\begin{equation}\label{jrproduct}
\binom{jr}{2}-j\binom{r}{2} = r\frac{jr(j-1)}{2}.
\end{equation}
Thus
\begin{alignat*}{2}
(\alpha(z)\beta(z))^{\binom{jr}{2}-  j\binom{r}{2}}
& = (\alpha(z)\beta(z))^{r \frac{jr(j-1)}{2}} &&
\mbox{ by \eqref{jrproduct} } \\
& =  ((\alpha(z)\beta(z))^r)^{\frac{jr(j-1)}{2}} && \\
& =  \beta(z)^{\frac{jr(j-1)}{2}} && \mbox{ by \eqref{24}. }
\end{alignat*}
Combining the previous equality with \eqref{firststep} and
\eqref{jrequation} 
we get
\begin{equation} 
\delta(x)^{j-1}= \beta(z)^{ \frac{jr(j-1)}{2}- \binom{j}{2}}=
\beta(z)^{\frac{j(j-1)(r-1)}{2}}.
\end{equation}

Since $j\neq1$, then $j-1 \not\equiv 0 \modu p$. So we may divide 
by $j-1$.  Therefore
\begin{equation}\label{lastone}
\delta(x)= \beta(z)^{\frac{j(r-1)}{2}}.
\end{equation}
Let $h \in \{1,\ldots, p-1\}$ such that $2h\equiv 1 \modu p$.
Then \eqref{lastone} implies \eqref{contradiction2}. 
\end{proof}

Suppose that 
Claim \ref{choosingcosetgood} (ii) holds. Then by Claim \ref{contradiction}
 we have that 
 $\delta(x)= \beta(z)^{hj(r-1)}$ and
$\delta(x) =\beta(z)^{hk(r-1)}$.
Thus $ hj(r-1)\equiv hk(r-1) \modu p$.
Since $r\not \equiv 1 \modu p$ and $2h\equiv 1\modu p$, we have that  
$k \equiv j\modu p $, a contradiction. Thus 
Claim \ref{choosingcosetgood}
(i) must hold. 

We apply now Claim \ref{contradiction} with  $\delta=1$. 
Thus $1= \delta(x) =  \beta(z)^{hj(r-1)}$.
Therefore $hj(r-1)\equiv 0 \modu p$. Since
$2h\equiv 1 \modu p$, either $j\equiv 0 \modu p$ or $r-1 \equiv 0 \modu p$. 
Neither is possible. That is our final contradiction.
\end{proof}
\end{section}

\begin{section}{Proof of Theorem A}

\begin{lemma}\label{countingabove}
 Let $G$ be a finite $p$-group and $N$ be a normal subgroup
of $G$. Assume that $\phi\in \Irr(N)$ is $G$-invariant. 
Then the set 
$\,\Irr(\, G\mid \phi\, )$ of all $\chi \in \Irr(G)$ 
lying over $\phi$ has either one 
or at     
least $p$ members. 
\end{lemma}
\begin{proof} Observe that if  $|G:N| = 1$, the claim holds 
 since $\Irr(\,G\mid \phi\, ) = {\phi}$ has one
member in that case. So we may assume that $|G:N| > 1$ and that the result
holds for all strictly smaller values of $|G:N|$.  
Since $G$ is a $p$-group,
there is some normal subgroup $M$ of $G$ such that $N < M$ and $|M:N| = p$.
Since $\phi$ is $M$-invariant and 
$M/N$ is cyclic, the set $ \Irr(\, M\mid \phi\, )$ consists of
exactly $p$ distinct extensions of $\phi$ to $M$.  There are two cases to
consider:
either some $\psi \in \Irr(\, M\mid \phi\, )$ is not $G$-invariant, or 
every $\psi \in \Irr(\, M\mid \phi\, )$ is $G$-invariant.

Assume that some $\psi \in \Irr(\, M\mid \phi\, )$ 
is not $G$-invariant.  Then the
stabilizer H of $\psi$ is a subgroup of index $p$ in $G$, and 
$\Irr(\, M\mid\phi\, )$
consists of the $p = |G:H|$ distinct $G$-conjugates of $\psi$.  
In this case
induction is a bijection of $\Irr(\, H\mid \psi\, )$
 onto $\Irr(\, G\mid \phi\, ) = \Irr(\, G\mid \psi\,)$.
Since $|H:M| = |G:N|/p^2 < |G:N|$, we know by induction that 
$\Irr(\,H\mid \psi\, )$
has either one or at least p elements. So the result holds in this case.

We may assume now that every $\psi \in \Irr(\, M\mid\phi\,)$ 
is $G$-invariant.  In this case
$\Irr(\, G\mid \phi\, )$ is the disjoint union of $p$ non-empty subsets
$\Irr(\,G\mid\psi\,)$,
for $\psi \in \Irr(\,M\mid \phi\,)$.  
Hence it has at last $p$ members, and the
result is proved.
\end{proof}

\begin{proof}[Proof of Theorem A]
Set $Z={\bf Z}(G)$. 
Assume that $2\eta(\chi,\psi) \leq p $. By Lemma \ref{main} we have that 
$\chi$ and $\psi$ vanish outside the center $Z$ of $G$.
Thus $\chi(1)\psi(1)=|G:Z|$ by
Corollary 2.30 of \cite{isaacs}. Let $\lambda\in \Irr(Z)$ and 
 $\mu\in \Irr(Z)$ be the unique characters lying below 
$\chi$ and $\psi$ respectively. Then $(\lambda \mu)$ lies below
$\chi\psi$ and
\begin{equation*} 
(\chi\psi)(g)= \left\{ \begin{array}{ll}
0 \mbox{ if $g \in G \setminus Z$}\\
|G:Z|(\lambda\mu)(g) \mbox{ if $g \in Z$.}
\end{array}
\right.
\end{equation*}
Thus $\chi\psi = (\lambda\mu)^G$. Since $2\eta(\chi,\psi)\leq p$, by 
Lemma  \ref{countingabove} applied to the normal 
subgroup $Z$ and the $G$-invariant 
character $\lambda\mu \in \Irr(Z)$,  we have that $\eta(\chi,\psi)=1$. 
\end{proof}
\end{section}

\begin{section}{Proof of Theorem B}
\begin{lemma}\label{center}
 Let $P$ be a finite  $p$-group, where $p$ is a prime number. 
 Let 
$\theta \in \Irr(P)$  be a character such      
that $\eta(\theta, \overline{\theta})=m>1$.
Then $p<2m+1$ and 
\begin{equation}\label{powersofp}
\theta(1) \in \{ p^i \mid i=1, \ldots, m-2\}
\end{equation}
\end{lemma}
\begin{proof} 
Since $G$ is supersolvable,
by  Theorem B of \cite{edith} we have  that  
$\theta(1)$ is the product of at most $\eta(\theta)-1$ prime numbers, where
$\eta(\theta)$ is the number of distinct non-principal
irreducible constituents of the product 
$\theta\overline{\theta}$. In our
case $\eta(\theta) = m-1$.  So $\theta(1)$ is the product 
of at most  $m-2$
prime numbers, each of which must be $p$. Hence \eqref{powersofp} 
holds.

Observe that $\Ker(\theta) = \Ker(\overline{\theta})$. Thus
 $\Ker(\theta{\overline{\theta}})\geq \Ker(\theta)\cap \Ker(\overline{\theta})=
\Ker(\theta)$. The group $P/\Ker(\theta)$ has a faithful character
$\chi$ such that the character 
$\theta\in \Irr(P)$ is inflated
from $\chi$. Also $\eta(\chi, \overline{\chi})=m$.
If $p$ is an odd prime,  by Theorem A we have that $\frac{p-1}{2}<m$. Thus
$p<2m+1$. 
\end{proof}

\begin{proof}[Proof of Theorem B]
Let
\begin{equation*}
S_n= \{ \, \prod (p_i)^{t_i} \mid  p_i< 2n+1,\, p_i 
\mbox{ is  a prime number for all
$i$, and } 0\leq t_i< n-2 \}
\end{equation*}
Let $p_1, p_2,\ldots, p_r$ be the set of distinct prime divisors of $|G|$. 
For $i=1, \ldots, r$, let $P_i$ be a Sylow $p$-subgroup of $G$.
Since $G$ is nilpotent, we have that 
$G = \prod_{i=1}^r P_i$. 
Thus the restriction
$\chi_{P_i}$ of $\chi$ to $P_i$ is
a multiple of an irreducible character $\theta_i$. Let
$m_i =\eta(\theta_i,\overline{\theta_i})$.
Observe that 
\begin{equation}\label{boundingmi}
m_i\leq n
\end{equation}
since 
$\chi\overline{\chi}= \prod_{i=1}^r \theta_i \overline{\theta_i}$
and $\eta(\chi,\overline{\chi})=n$. Also
\begin{equation}\label{gradochi}
\chi(1) = \prod_{i=1}^r \theta_i(1).
\end{equation}
If $m_i=1$ then the character $\theta_i$ is linear. 
Assume that $m_i>1$. By Lemma  \ref{center} we have that 
$\theta_i(1)= p^{t_i}$, where  
$t_i\leq m_i-2$, and  $p_i< 2m_i+ 1$. By \eqref{boundingmi} we have 
$p_i< 2n+1$. So by \eqref{gradochi}
we have
\begin{equation*}
\chi(1)\in S_n.
\end{equation*}
\end{proof}

\end{section}

\begin{section}{Examples}

\begin{proposition} 
Let $p$ be an odd prime and $m$ be any positive integer.
There exist a $p$-group $G$ and a faithful character  $\chi\in \Irr(G)$
such that  $\chi\chi$ has $\frac{p+1}{2}$ distinct irreducible
constituents, $\chi$ does not vanish on $G \setminus {\bf Z}(G)$
and $\chi(1)=p^m$.
\end{proposition}
\begin{proof} 
Let $C =\{0,1,2,\ldots, p-1\}$
be the additive group of integers modulo $p$. Let $E$ be an extraspecial 
of exponent $p$ and order $p^{2m-1}$ if $m>1$. Otherwise let
$E$ be a cyclic group  of order $p$. 
  
Set 
$$A=E^{C}= \mbox{ the set of all functions from $C$ to $E$}.$$
Observe that $A$ is just the direct product of $p$ copies of $E$.
and has order $p^{p(2m-1)}$.  Also 
observe that $C$ acts on $C$ with the regular action,
and thus it acts on the group $A$ with action
	$$(f^c)(x)= f(c+x)$$ 
\noindent for any $c, x \in C$. 

Let $G$ be the semi-direct product of 
$C$ and $A$. Observe that $|G|= p^{p(2m-1)+1}$ and ${\bf Z}(G)$ is cyclic 
of order $p$. Thus $|G:{\bf Z}(G)|=p^{p(2m-1)}$. Since $p$ is odd, 
$p^{p(2m-1)}$ is not a square. So no $\chi \in \Irr(G)$ can vanish 
on $G\setminus {\bf Z}(G)$ (see Corollary 2.30 of \cite{isaacs}). 

Let $\lambda \in \Irr(E)^{\#}$ have degree $p^{m-1}$.
If $m>1$, such character exists
 since $E$ is an extra-special group of order $p^{2m-1}$. 
Observe that 
$\lambda\lambda= \lambda(1) \mu$ for some 
 $\mu\in \Irr(E)$. For each $i=0,1,\ldots p-1$, let
 $\mu_i \in \Irr(A)$ be a character such that
for all $a\in A$ we have 
$$\mu_i(a)= \mu(a(i)).$$
For $i=0, 1, \ldots, p-1$, choose  $\lambda_i \in \Irr(A)$ so that we have,
for all $a \in A$,
	$$\lambda_i(a)= \lambda (a(i)).$$
Observe that the stabilizer of $\lambda_i $ is $A$. 
Thus $\lambda_i^G \in \Irr(G)$.  Set
$\chi = \lambda_0^G \in \Irr(G)$.

\begin{claim}\label{counting} Let $S= \{1, 2, \ldots, \frac{p-1}{2}\}$. Then 

(i) The set $\{\lambda_0 \lambda_i \mid i\in S\}\cup \{\mu_0\}$
is a subset of the irreducible constituents of $(\chi \chi)_A$. 

(ii) If $i, j \in S$ and $i \neq j$, then 
$\lambda_0 \lambda_i$ and $\lambda_0 \lambda_j$
 are not $G$-conjugate. 

(iii) Given 
 $k, l=0,1,2, \ldots, p-1$, where $k \neq l$, there exist
$j \in S$ such that $\lambda_k\lambda_l$ and $\lambda_0\lambda_j$
are $G$-conjugate. Also for any $k \in S$, 
$\lambda_k\lambda_k$ is $G$-conjugate to 
$\lambda_0\lambda_0$. 
\end{claim}
\begin{proof}
We can check that 
\begin{equation}\label{productchi}
\chi_A= \sum_{i=0}^{p-1} \lambda_i.
\end{equation}
Thus
\begin{equation*}
(\chi\chi)_A= \chi_A\chi_A= \sum_{i,j=0}^{p-1} \lambda_i\lambda_j.
\end{equation*}
By definition of $A$, we have that  $\lambda_0 \lambda_i \in \Irr(A)$ for 
any $i=1, \ldots, p-1$. Therefore 
$\lambda_0\lambda_i$ is  an irreducible constituent
of $(\chi \chi)_A$  for any $i \in S$. Also $\mu_0$ is an irreducible
constituent of $\chi\chi$ since $\lambda_0\lambda_0 =\lambda_0(1)\mu_0$. 

Assume that $\lambda_0 \lambda_i$ and $\lambda_0 \lambda_j$, where
$i, j \in S$ and $i\neq j $,  are 
$G$-conjugates. Then there exists $c\in C$ such that 
$c+0\equiv j\modu p$ and $c+i\equiv 0\modu p$. Thus $i+j\equiv 0\modu p$.
Therefore either $i>\frac{p-1}{2}$ or $j>\frac{p-1}{2}$. 
We conclude that either $i \not\in S$ or $j \not\in S$.

Fix $k, l =0, 1, \ldots, p-1$, where $k \neq l$. 
Let $c\equiv -k \modu p$ and 
$r\equiv l-k \modu p$, where
$c, r \in \{0,1,2, \ldots, p-1\}$. Observe that $r\not\equiv 0\modu p$
since $k \neq l$.  Then 
$(\lambda_k \lambda_l)^c= \lambda_0 \lambda_r$. If $r \in S$, we  have 
finished. Otherwise
 $r>\frac{p-1}{2}$. Let $t\equiv  -r \modu p$, where
$t\in  \{0,1,2, \ldots, p-1\}$. Observe that
$t\leq \frac{p-1}{2}$ and $(\lambda_0 \lambda_r)^t= \lambda_0\lambda_t$.

	Clearly $(\lambda_0 \lambda_0)^c= \lambda_c \lambda_c$  for
any $c\in C$, and thus
the last statement of the claim follows. 
\end{proof}

Observe that the stabilizer of $\lambda_0\lambda_i$, where $i \in S$, is
A. Thus $(\lambda_0\lambda_i)^G \in \Irr(G)$. 
By Claim \ref{counting}, we have that 
$(\lambda_0\lambda_i)^G \neq (\lambda_0\lambda_j)^G$ if 
$i \neq j$ and $i, j \in S$. Also
$$[(\lambda_0\lambda_i)^G, \chi\chi]= [(\lambda_0\lambda_i), 
(\chi\chi)_A]\neq 0, $$
\noindent where the last inequality follows from \eqref{productchi}.
Similarly, we can check that $(\mu_0)^G \in \Irr(G)$.
It follows from Claim \ref{counting} that 
$\{(\lambda_0\lambda_i)^G\mid i \in S\}\cup \{(\mu_0)^G\}$ 
is the set of all distinct 
 irreducible
constituents of $\chi\chi$. This 
 set  has $|S|+1= \frac{p-1}{2}+1= \frac{p+1}{2}$
 elements.
\end{proof}

\end{section}

{\bf Acknowledgment.} Professor Everett C. Dade 
made a substantial contribution to this work. 
I thank him for that.

\bibliographystyle{amsplain}

\begin{thebibliography}{9}


\bibitem{edith} E. Adan-Bante, Products of characters and derived length, 
to appear in J. Algebra. 

\bibitem{edith2} E. Adan-Bante, Products of characters with few irreducible
constituents, preprint. 

\bibitem{us} E. Adan-Bante, M. Loukaki, A. Moreto, 
Homogeneous products of characters, 
preprint. 



\bibitem{isaacs} I.M.Isaacs, Character Theory of Finite Groups. New York-San
Francisco--London: Academic Press 1976


\end{thebibliography}

\end{document}